\documentclass[11pt]{amsart}

\usepackage[left=3cm,tmargin=3cm,centering]{geometry}

\usepackage{xcolor}
\usepackage[backref=page]{hyperref}
\definecolor{couleur_cite}{rgb}{0.05,.4,0.05}
\definecolor{couleur_link}{rgb}{0.05,0.05,0.4}
\hypersetup{
bookmarksopen,bookmarksnumbered,
colorlinks,
linkcolor=couleur_link,citecolor=couleur_cite,
}

\usepackage[initials]{amsrefs}  
\usepackage{bm,amsmath,amssymb,amsthm,amscd}

\usepackage[all]{xy}
\usepackage{stmaryrd}

\usepackage{amsfonts}

\usepackage{abbreviations}

\newcommand{\sumf}{\sum_{f\in \CmB_2(N)}}

\newcommand{\sumc}{\sum_{c\equiv 0 (N)}}
\newcommand{\trd}{\MTr' \Delta^{-1}}

\usepackage[latin1]{inputenc}

\numberwithin{equation}{section}
\usepackage{longtable}    
\allowdisplaybreaks[1]    

\theoremstyle{plain}
\newtheorem{theorem}{Theorem}

\newtheorem{proposition}{Proposition}[section]
\newtheorem*{proposition*}{Proposition}
\newtheorem{lemma}{Lemma}[section]
\newtheorem*{lemma*}{Lemma}
\newtheorem*{corollary*}{Corollary}

\theoremstyle{definition}

\newtheorem{remark}{Remark}

\theoremstyle{remark}

\numberwithin{equation}{section}

\begin{document}

\title[]{Remarks on a special value of the Selberg zeta function} 

\author[]{Nicolas Templier}   
\address{Institute for Advanced Study, School of Mathematics, 08540 Princeton, NJ, USA}
\email{nicolas.templier@normalesup.org}

\date{\today}

\begin{abstract}
Let $\CmZ_{Y_0(N)}$ be the constant term of the logarithmic derivative at $s=1$ of the Selberg zeta function of the modular curve $Y_0(N)$. Jorgenson and Kramer established the bound $\CmZ_{Y_0(N)}=O_\epsilon(N^\epsilon)$, $\epsilon>0$ by relating it to geometric invariants. In this article we give, for $N$ prime, another proof via $L$-functions and exponential sums improving on a previous approach by Abbes-Ullmo and Michel-Ullmo. We further derive a power of $\log N$ bound along the same line.
\end{abstract}

\maketitle


\tableofcontents


\section{Introduction}
\subsection{Selberg zeta function.}
Let $X$ be the quotient of the upper half-plane by a fuchsian group of the first kind, equipped with its hyperbolic metric. Let $Z_X(s)$ be the Selberg zeta function associated to $X$, see~\cite{book:Iwan02}*{\S 10.8}. It has a meromorphic continuation to $\BmC$ and a functional equation, as follows from the Selberg trace formula. We are interested in the quantity:
\begin{equation}\label{def:zx}
 \CmZ_X:=\lim_{s\rightarrow 1} \frac{Z'_X(s)}{Z_X(s)}-\frac{1}{s-1}=\frac{1}{2}\frac{Z''_X}{Z'_X}(1).
\end{equation}
We give briefly in the last section~\ref{sec:rem} some reasons why to study this quantity. It appears explicitly in the self-intersection of the dualizing sheaf (Abbes-Ullmo~\cite{AU97}), as a regularized trace of $\Delta^{-1}$ (Steiner~\cite{Stei05}), in the study of Faltings delta function (Jorgenson-Kramer~\cite{JK06}) and in the asymptotic height of Heegner points (\cite{Temp:shifted}). The \emph{arithmetic significance} of $\CmZ_X$ remains mysterious, at least for the author. We believe it deserves further study which was a main motivation for the present paper.

Here we shall consider the particular case $X=Y_0(N)$, the modular curve and the behavior as $N\to \infty$. In~\cites{Jorg01}, Jorgenson and Kramer proved that for $N$ squarefree:
\begin{equation}\label{JK}
 \CmZ_{Y_0(N)}=O_\epsilon(N^{\epsilon}).
\end{equation}
Their method was geometric, relating $\CmZ_X$ to several hyperbolic invariants of $X$, see \S~\ref{sec:remarks:JK} below for a summary. A variant works also for the study of $\log Z'_X(1)$, the Faltings delta function $\delta\subt{Fal}(X)$, for the compactified $X_0(N)$, and for $N$ not necessarily squarefree, see~\cite{JK06}.

Previously, Michel and Ullmo~\cite{MU98} established the weaker\footnote{but stronger than a \Lquote{na\"ive} bound which would be $O_\epsilon(N^{1+\epsilon})$ and the result is good enough for the application to asymptotics of the self-intersection of the dualizing sheaf} bound $\CmZ_{Y_0(N)}=O_\epsilon(N^{7/8+\epsilon})$ also for $N$ squarefree. Their method is analytic and relies on the theory of Rankin-Selberg $L$-functions.

It seems interesting to compare the two results, a natural question is as follows. \emph{Is it possible to compare the two methods of proof and explain the disparity of exponents?} In this article we answer this question by showing that, for $N$ prime, an improvement of the argument of Michel and Ullmo yields the same bound~\eqref{JK}. Hence the \Lquote{geometric} and \Lquote{analytic} approaches truly yield bounds of comparable strength for $\CmZ_{Y_0(N)}$.

\subsection{Main result.}
For simplicity we have assumed $N$ prime. It seems the argument should work for arbitrary $N$, but this would demand technical extra care at many steps, typically in the definition of the symmetric square $L$-function and in the combinatoric of newforms (however see~\cite{ILS00} for a good treatment of this kind of questions for $N$ square-free). We want to avoid these complications here.

Let $\CmB_2(N)$ be an orthogonal basis of forms of weight $2$, level $N$ and trivial Nebentypus. We demand that the forms $f\in\CmB_2(N)$ are Hecke eigenforms and are normalized with first Fourier coefficient $1$ (since $N$ is prime these forms are new). Let $g(N)=\# \CmB_2(N)$ be the genus. The aim of this article is to prove:
\begin{theorem}\label{th:main} For $N$ prime, and $\epsilon>0$ one has:
\begin{equation}
  \frac{1}{g(N)}\sum_{f\in \CmB_2(N)}
 \frac{L'}{L}(1,\MSym^2 f)
 =2\frac{\zeta'}{\zeta}(2)
 +O_\epsilon(\frac{1}{N^{1-\epsilon}}).
\end{equation}
\end{theorem}
\remark The bound $O_\epsilon(N^{-1+\epsilon})$ improves greatly on the $O_\epsilon(N^{-1/8+\epsilon})$ by Michel-Ullmo~\cite{MU98}, which they obtained in a slightly different way, see \S~\ref{sec:rem:MU} below for a comparison. An important difference is that we apply Poisson summation formula (\S \ref{sec:proof:poisson}) instead of Weil's bound for Kloosterman sums.

\remark In Theorem~\ref{th:log} from \S~\ref{sec:rem:log} we improve the remainder term into $O(\frac{(\log N)^7}{N})$. Since the arguments are technical and not really illuminating we decided to postpone this improvement to the last section of the paper. Together with the identity~\eqref{AU} below this shows that for $N$ prime:
\begin{equation}\label{Zlog}
 \CmZ_{Y_0(N)}\ll (\log N)^7.
\end{equation}
It should be noted that~\cite{Jorg01} already derived the lower bound $\CmZ_{Y_0(N)} \ge O(\log N)$ (it is not known whether $\CmZ_{Y_0(N)}$ is positive or negative).

\remark It would be interesting to obtain an asymptotic for $\CmZ_{Y_0(N)}$. However this task seems difficult in view of the number of sums and integrals that are performed. It seems to the author more ideas need to be introduced. Especially unbinding the $l$ and $c$-sums is a delicate matter. Probably it would be necessary to treat the archimedian integrals in a different way. A first step would be to insert~\eqref{def:Vs} inside the Fourier-Bessel integral~\eqref{def:g} of $g(l;c)$ and derive sharper estimates, in some sense an archimedian analogue of the evaluation of the complete exponential sum $R(l,n;c)$ in \S~\ref{sec:proof:exp}. 

\subsection{An identity.} Cutting a long story short, from Abbes-Ullmo \cite{AU97}*{Proposition 3.3.6 p. 58} (see also \cite{MU98}*{eq. (11)}) one has for $N$ prime:
\begin{equation}\label{AU}
 \frac{-1}{2g}\CmZ_{Y_0(N)}=
 \frac{1}{g(N)}\sum_{f\in \CmB_2(N)}
 \frac{L'}{L}(1,\MSym^2 f)
 -2\frac{\zeta'}{\zeta}(2)
 +O(\frac{\log N}{N}).
\end{equation}

\remark Let us observe that the remainder term is in fact explicit in~\cite{AU97} at least when $N\equiv 11\qmod{12}$, so that we truly have an \emph{identity}. When $N$ is only square-free the $\log N$ needs to be modified by several divisor-like functions. It is perhaps possible that these terms might be interpreted in terms of $L$-functions and newforms, and would match with the variant of Theorem~\ref{th:main}. It would be interesting to know whether $\CmZ_{Y_0(N)}$ is still bounded by a power of $\log N$ in the general squarefree case.

\remark The fact that the geometric quantity $\CmZ_{Y_0(N)}$ is related to a quantity which involves weight two forms sounds a bit strange. There is no particular reason to believe that the bundle $\Omega^1_{Y_0(N)}$ would play a central role in this riemannian setting. The quantity $\CmZ_{Y_0(N)}$ is related to a regularized sum of eigenvalues of Maass wave forms and is a spectral invariant of the metric, see~\S~\ref{sec:rem:spectral}. It would be interesting to clarify the proof of~\eqref{AU} in~\cite{AU97} and try to make use of Maass forms instead of modular forms.

\remark The use of $L$-functions is an intermediate object for the bound~\eqref{Zlog}. If one inspects the proof given below, we really have an identity (sophisticated) between $\CmZ_{Y_0(N)}$ and a large combination of exponential sums $R(l,n,c)$. It would be interesting to give an intrinsic proof of such an identity.

\subsection{Brief outline of proof.} The first step is to express $\frac{L'}{L}(1,\MSym^2 f)$ in a convenient fashion, which we do by exploiting the functional equation of the symmetric square $L$-function. Although short, the argument is tricky since we need to move the integration path \emph{twice}.
It is more usual (and natural) to apply the approximation functional equation method to values at $1/2$ (center of the critical strip) rather than at $1$ (edge of the critical strip). Indeed at $1$, it is usually the logarithm of the conductor that really matters and many of the deep and delicate questions are raised and solved in that scale (e.g. the prime number theorem). For the aim of this article where we seek for ample cancellations, this truncation at the square-root of the conductor is well-suited.

The average $\sumf$ is dealt with Petersson formula (\S \ref{sec:proof:petersson}), this is very classical. Then instead of applying Weil's bound for Kloosterman sums, we make use of the (relatively) large length of the $n$-sum. Applying Poisson summation formula, we create certain \emph{complete} exponential sums of dimension $2$. 

It is easy to prove square-root cancellation for these sums (\S \ref{sec:proof:exp}). Because the dimension is larger, the saving is better than for Kloosterman sums. This enables to conclude the proof of the Theorem.

\remark Formally several steps in our proof are similar to the article~\cite{IM01} by Iwaniec-Michel, where a sharp bound for the average of $L(1/2,\MSym^2 f)^2$ is provided. In particular the exponential sum $R(l,n;c)$ appears also in~\cite{IM01} and its explicit value is exploited there to gain cancellations via Heath-Brown quadratic large sieve inequality. In the present paper the parameters get restricted to a smaller range and a square-root cancellation of the exponential sum is sufficient for our purpose.

\subsection{Notation.}\label{sec:intro:notation} The index of the congruence subgroup $\Gamma_0(N)$ in $\Gamma(1)=SL_2(\BmZ)$ is $N+1$. One has $-\frac{\zeta'_N}{\zeta_N}(s)=\frac{\log N}{N^s-1}$ and we shall use the fact that $g(N)\gg N$ (of course the exact value of $g(N)$ is known by Riemann-Roch). From line to line the arbitrary $\epsilon>0$ may change its value. The Petersson inner product is:
\begin{equation}
 (f,g)=\int_{Y_0(N)} f\overline{g} dxdy.
\end{equation}

The symmetric square of $f$ is denoted $\MSym^2 f$. Recall that its $L$-function is of degree $3$, conductor $N^2$, root number $+1$, that the $L$-series is:
\begin{equation}\label{def:sym2}
 L(s,\MSym^2 f)=\zeta^{(N)}(2s)\sum^\infty_{n=1} \frac{\lambda_f(n^2)}{n^s},
\end{equation}
and its gamma factor is $\gamma(s):=\Gamma_\BmR(s+1)^2\Gamma_\BmR(s+2)$, where $\Gamma_\BmR(s)=\pi^{-s/2}\Gamma(\frac{s}{2})$.

One has: $L(1,\MSym^2 f)=\frac{8\pi^2}{N}(f,f)$ and $g(N)=\frac{N+1}{12}+O(1)$. We set:
\begin{equation}
 \omega_f:=L(1,\MSym^2 f)^{-1}.
\end{equation}
From~\cite{HL94} one has $N^{-\epsilon}\ll_\epsilon \omega_f\ll_\epsilon N^{\epsilon}$, but we won't need this fact.

As usual $e(y):=e^{2i\pi y}$ and $e_c(y):=e(\frac{y}{c})$. The Kloosterman sum is:
\begin{equation}
 S(m,n;c):=\sum_{x\in (\BmZ/c\BmZ)^{\times}}
 e_c(m\overline{x}+nx).
\end{equation}
We shall sometimes abbreviate $x\in \BmZ/c\BmZ$ into $x(c)$ and $x\in (\BmZ/c\BmZ)^{\times}$ into $x(c)^*$.

We shall not use the notation $C_F$ from~\cites{AU97,MU98} which plays the role of the right-hand side of~\eqref{AU}. Thus for the convenience of the reader we briefly explain how to express $C_F$ in our notation. The quantity $\vol\cdot g\cdot C_F$ was defined in~\cites{AU97,MU98} to be the constant term in the Laurent expansion at $s=1$ of:
\begin{equation}
 \vol \int_{Y_0(N)} E_{\infty}(s,z) \sum_{f\in \CmB_2(N)} \frac{\abs{f(x+iy)}^2}{(f,f)} dxdy.
\end{equation}
A standard unfolding computation then yields:
\begin{equation}
 \vol\cdot g\cdot C_F=\sum_{f\in \CmB_2(N)}
 \frac{L'}{L}(1,\MSym^2 f) + g[-2\frac{\zeta'}{\zeta}(2)+1-\log 4\pi +\frac{\log N}{N+1}].
\end{equation}

\section{Proof}\label{sec:proof}
\subsection{Approximate functional equation at the edge.}\label{sec:proof:approximate} From the analytic properties of $L(s,\MSym^2 f)$ previously recalled in \S \ref{sec:intro:notation}, one may deduce:
\begin{proposition}\label{lem:approximate} Fix $u\mapsto G(u)$ a meromorphic, even function on $\BmC$, with only pole at $u=0$ of Laurent expansion $\frac{1}{u^2}+O(1)$, and bounded on the vertical strip $-4 \le \MRe u\le 4$, $\abs{\MIm u}\ge 1$. Then, for all $f\in\CmB_2(N)$:
\begin{multline}
\frac{L'}{L}(1,\MSym^2 f)=
\frac{I^*_f(1)-I^*_f(0)}{\gamma(1)L(1,\MSym^2 f)} 
-\log N -\frac{\gamma'}{\gamma}(1)+\\
+\frac{\zeta(2)}{L(1,\MSym^2 f)}
\biggl(\frac{\zeta'}{\zeta}(2)
+\frac{\gamma'}{\gamma}(1)
+\log N+O(\frac{\log N}{N})
\biggr) 
\end{multline}
where $I^*_f(s)$ is the following series:
\begin{equation}\label{def:Ifs}
 I^*_f(s):=\frac{1}{N}\sum^{\infty}
 _{\substack{m=1\\ (m,N)=1}}
 \sum^\infty_{n=2}
 \lambda_f(n^2)
 (\frac{N}{m^2n})^s
 V_s(\frac{m^2n}{N})
\end{equation}
and where $V_1(y)$ and $V_0(y)$ are the smooth functions defined by:
\begin{equation}\label{def:Vs}
 V_s(y):=\int_{(3)} y^{-u} \gamma(s+u)G(u)\frac{du}{2i\pi}.
\end{equation}
\end{proposition}
\begin{proof} Put $\Lambda(s):=N^s\gamma(s)L(s,\MSym^2 f)$ and consider the integral $I(1)$ where:
\begin{equation}
 I(s):=\int_{(3)}
 \Lambda(s+u) G(u)
 \frac{du}{2i\pi}.
\end{equation}
We may move the integration contour onto $\Re u=-3$ crossing a pole at $u=0$ of residue $\Lambda'(1)$, and apply the functional equation $\Lambda(s)=\Lambda(1-s)$:
\begin{equation}
 I(1)=\Lambda'(1) + I(0).
\end{equation}
It is clear that:
\begin{equation}
 \Lambda'(1)=N\gamma(1)L(1,\MSym^2 f) 
 \Bigl(\log N +\frac{\gamma'}{\gamma}(1)+\frac{L'}{L}(1,\MSym^2 f)
 \Bigr).
\end{equation}

Expanding $L(s,\MSym^2 f)$ into series, we get:
\begin{equation}
 I(s)=\sum^{\infty}
 _{\substack{m=1\\ (m,N)=1}}
 \sum^\infty_{n=1}
 \lambda_f(n^2)
 (\frac{N}{m^2n})^s
 V_s(\frac{m^2n}{N}),
\end{equation}
from which we deduce:
\begin{equation}
 I(s)=N I^*_f(s)+\int_{(3)} N^{s+u}\zeta^{(N)}(2s+2u)\gamma(s+u)G(u) \frac{du}{2i\pi}.
\end{equation}
Observe that $N^{s+u}\zeta^{(N)}(2s+2u)=(N^{s+u}-N^{-s-u})\zeta(2s+2u)$ is never too small as $N\to \infty$ (at most bounded). For $s=0$, we move back the integration to $\MRe(u)=\epsilon$ crossing a pole at $u=\Mdemi$ and getting:
\begin{equation*}
 I(0)=NI^*_f(0)+O_\epsilon(N^{\epsilon}).
\end{equation*}
It one instead circles the double pole at $u=0$ at distance $1/\log N$ it is possible to improve the latter into:
\begin{equation}\label{I-Istar}
 I(0)=NI^*_f(0)+O(\log N).
\end{equation}

For $s=1$, we move the integration to $\MRe(u)=-1$, crossing a pole at $u=0$, getting:
\begin{equation}
 I(1)=N I^*_f(1)+(N-\frac{1}{N})\zeta(2)\gamma(1)
 \Bigl(\frac{\zeta'}{\zeta}(2)+\frac{\gamma'}{\gamma}(1)+ 
 \frac{N+N^{-1}}{N-N^{-1}}
 \log N
 \Bigr)
 +O(1).
\end{equation}
\end{proof}

A straightforward consequence of Petersson formula (recalled in the next section) and Weil's bound for Kloosterman sums is the following (see, e.g,\cite{MU98}*{eq. (15)}):
\begin{lemma}\label{lem:omega} The weight $\omega_f=L(1,\MSym^2 f)^{-1}$ satisfies:
 \begin{equation}
  \frac{1}{g(N)} \sumf \omega_f
  =\frac{6}{\pi^2}+O(\tau(N)N^{-3/2}).
 \end{equation}
 \qed
\end{lemma}

In view of Proposition~\ref{lem:approximate} and Lemma~\ref{lem:omega}, we need to estimate, for $s=0,1$:
\begin{equation}
 \frac{A_s}{g(N)}:=
 \frac{1}{g(N)}\sumf
 \omega_f I^*_f(s)
 ,\quad
 \text{as $N\to \infty$}
\end{equation}
which is the aim of the following paragraphs.

\remark The presence of the $m$-sum in~\eqref{def:Ifs} does not play any important role in the sequel and the reader may consider $m=1$ as the typical case. The author believes it should be possible in this section to treat in a more elegant way the term $\zeta^{(N)}(2s)$ in~\eqref{def:sym2} -- from which the $m$-sum arises. However it is not clear how one should modify the approximate functional equation method.

\subsection{Applying Petersson formula.}\label{sec:proof:petersson}
For $m,n\ge 1$ let:
\begin{equation}
 \Delta^*_N(m,n):=
 (4\pi)^{-1}
 \sum_{f\in \CmB_2(N)}
 \frac{\lambda_f(m)\lambda_f(n)}{(f,f)}
\end{equation}
and
\begin{equation}
 \Delta_N(m,n):=\frac{1}{g(N)}
 \sum_{f\in \CmB_2(N)} \omega_f
 \lambda_f(m)\lambda_f(n).
\end{equation}
Since $N$ is prime, the relation between $\Delta_N$ and $\Delta^*_N$ is rather trivial\footnote{for the general case which involves combinatorics on newforms, see~\cite{ILS00}}, see \S~\ref{sec:intro:notation}:
\begin{equation}
 \Delta_N(m,n)=(\frac{1}{2\pi^2}+O(\frac{1}{N}))
 \cdot
 \Delta^*_N(m,n).
\end{equation}

The Petersson formula in weight $2$ reads~\cite{book:IK04}*{\S~14.10}:
\begin{equation}
 \Delta^*_N(m,n)=\delta_{m,n}-2\pi \sum_{c\equiv 0 (N)} \frac{S(m,n;c)}{c} J_1(4\pi\frac{\sqrt{mn}}{c}).
\end{equation}
The $c$-sum converges absolutely by the Weil estimate for Kloosterman sums (or any power saving improvement on the trivial bound) and the bound $J_1(x)\ll \min(x,x^{-1/2})$ (see e.g.,~\cite{book:Iwan02}*{Appendix B}).

It thus remains to estimate ($s=0,1$):
\begin{equation}
 A_s=\frac{1}{N}\sum^{\infty}
 _{\substack{m=1\\ (m,N)=1}}
 \sum^\infty_{n=2}
 (\frac{N}{m^2n})^s
 V_s(\frac{m^2n}{N})
 \sumc
 \frac{S(n^2,1;c)}{c} J_1(\frac{4\pi n}{c}).
\end{equation}
\subsection{Truncation of the $c$-sum.}
Before proceeding further it is important to tail the sum over the integers $c$. Choose $C:=N^{10}$, then it is not difficult to see that the $c> C$ yields a negligible contribution to $A_s$ (certainly much smaller that $O(\frac{1}{N})$). From now on the $c$-sum will be always tacitly restricted to $c\le C$.

\subsection{Applying Poisson formula.}\label{sec:proof:poisson} Let $W_s(y):=y^{-s}V_s(y)$ and $X:=N/m^2$. Exchanging the $n$ and the $c$-summation, we ought to estimate:
\begin{equation}\label{def:Bs}
 B_s(X):= \frac{1}{X}
 \sum_{c\equiv 0 (N)}
 \sum^\infty_{n=2} 
 \frac{S(n^2,1;c)}{c} W_s(\frac{n}{X})J_1(\frac{4\pi n}{c}).
\end{equation}
As in~\cite{IM01}*{\S 3} we may replace, up to logarithmic factors, the truncation $W_s(\frac{n}{X})$ by $U(\frac{n}{Y})$, where $U\in \CmC^\infty(1,2)$, $1 \le Y\le XN^{\epsilon}$ and $U^{(i)}\ll_i 1$ for all $i\in\BmN$. We denote the resulting sum by $B(Y)$.

We apply Poisson summation formula (see, e.g, \cite{book:IK04}*{\S 4.3}) to the $n$-sum, filling in the residue classes modulo $c$. We get:
\begin{equation}\label{Bs-complete}
 B(Y)=\frac{1}{Y}
 \sum_{l\in \BmZ}
 \sum_{c\equiv 0 (N)}
 \frac{R(n,l;c)}{c}
 g(l,c)
\end{equation}
where we have introduced the following complete exponential sum:
\begin{equation}
 R(n,l;c):=\sum_{n\qmod{c}}
 S(n^2,1;c)e_c(ln)
\end{equation}
and the Fourier transform:
\begin{equation}\label{def:g}
 g(l,c):=\int^\infty_{-\infty}
 e(-\frac{lx}{c})U(\frac{x}{Y})J_1(\frac{4\pi x}{c}) \frac{dx}{c}
\end{equation}
(note that since $U$ is supported on $(1,2)$, the $J_1$-function is evaluated at a positive number).

\subsection{Truncation of the spectral sum.}\label{sec:proof:truncation}
It is clear that $g(l,c)\ll Y^{2}/c^2$ and integrating by parts several times, we get more precisely ($s=0,1$):
\begin{equation}\label{glc-decay-l}
 g(l,c)\ll_{i} \frac{Y^2}{c^2} 
 (1+\frac{c}{Y})^{i}
 l^{-i} ,
 \qtext{for all $i\in \BmN,\ l\not=0$}.
\end{equation}
Thus, for any $\epsilon>0$, the terms in~\eqref{Bs-complete} with $\abs{l}> \dfrac{c}{Y}
N^{\epsilon}$ have a negligible contribution, by taking $i$ large enough (depending on $\epsilon$ only).

\subsection{Square-root cancellation.}\label{sec:proof:exp} By definition of the Kloosterman sum, one  has:
\begin{equation}\label{def:R}
 R(l,n;c)=\sum_{n(c)}\sum_{x(c)^*}
 e_c(n^2\overline{x}+x+ln).
\end{equation}
Making the (bijective) change of variable $(n,x)\leadsto (xn,x)$, we get:
\begin{equation}
 R(l,n;c)=\sum_{n(c)}
 \sum_{x(c)^*}
 e_c(x(n^2+ln+1))
\end{equation}
in which we recognize a Ramanujan sum as the inner summation. We thus may write:
\begin{multline}\label{squareroot}
 R(l,n;c)=\sum_{n(c)}
 \sum_{a|c} a\mu(\frac{c}{a})\delta(a|n^2+ln+1)
 \\
 = c
 \sum_{a|c} \mu(\frac{c}{a})
 \times
 \#\{n(a),\ n^2+ln+1\equiv 0(a)\}\ll c\tau(c)^2 \ll_\epsilon c^{1+\epsilon}.
\end{multline}

\subsection{Conclusion.}\label{sec:proof:conclusion} Returning to~\eqref{Bs-complete}, we get from~\S \ref{sec:proof:truncation} and~\eqref{squareroot}:
\begin{equation}\label{Bs}
 B(Y)\ll_\epsilon \frac{1}{Y}
 \sum_{c\equiv 0(N)} \frac{c}{Y}N^\epsilon \times c^\epsilon \times \frac{Y^2}{c^2}
 \ll_\epsilon \frac{1}{N} C^\epsilon N^{2\epsilon}.
\end{equation}
The same is true for $B_s(X)$, thus we conclude:
\begin{equation}
 A_s=\sum_{(m,N)=1} \frac{1}{m^2} B_s(N/m^2)= O_\epsilon(N^\epsilon).
\end{equation}

\section{Further remarks}\label{sec:rem}
\subsection{Another approximation of the $L$-values.}\label{sec:rem:MU}
We review in detail the main steps in Michel-Ullmo~\cite{MU98} since we believe the comparison is instructive. We simplify somehow the exposition in order to discuss the main points and the reader should refer to~\cite{MU98} for precise statements.

Recall the relation between the symmetric square and the Rankin-Selberg convolution (which follows from Hecke relations):
\begin{equation}
 L(s,\MSym^2 f)\zeta^{(N)}(s)=L(s,f\times f)=\zeta^{(N)}(2s)
 \sum^\infty_{n=1} \frac{\lambda_f(n)^2}{n^s}.
\end{equation}
Thus it is natural to consider the function:
\begin{equation}
 L_f(s):=L(s,f\times f)-\zeta^{(N)}(s)L(1,\MSym^2 f)
\end{equation}
which has a holomorphic continuation to all of $\BmC$. The quantity one needs to estimate is then:
\begin{equation}
 \frac{1}{g(N)} \sumf
 \frac{L'}{L} (1,\MSym^2 f) \approx
 \frac{1}{g(N)} \sumf
 \omega_f L_f(1).
\end{equation}

It is not easy to express $L_f(1)$ as a truncated sum of eigenvalues, so~\cite{MU98} proceed in the following way. They consider $L_f(\sigma)$, or more generally $L_f(s)$ for $\MRe s=\sigma$ fixed and near $1$. Then one applies the Phragm\'en-Lindel\"of principle to reach the value at $1$. Here we use the approximate functional equation method instead (\S~ \ref{sec:proof:approximate}). 

The next step is similar to our \S~\ref{sec:proof:petersson}: one applies Petersson formula to the average of $L_f(\sigma)$. The formula they obtain is quite complicated because there are several truncations (which was the advantage of our more flexible Proposition~\ref{lem:approximate}). The quantity \Lquote{$T'_1$} in \cite{MU98}*{\S~3.1} should be considered as the heart matter. One may observe that it is not very far from our quantity $B_s(X)$ defined in~\eqref{def:Bs}, except that the Kloosterman sum is \Lquote{$S(n,n;c)$} instead of $S(n^2,1,c)$.

The next step is different: we apply Poisson summation formula while~\cite{MU98} applies Weil's bound for Kloosterman sums. As explained in the introduction, this yields much less savings. Taking also into account the remaining terms which are difficult to estimate, this explains why the exponent $7/8$ obtained in~\cite{MU98} is really much larger.

The exponential sums that would arise from~\cite{MU98} by applying Poisson are the following. It clearly also has a square-root cancellation:
\begin{equation}\label{exp-MU}
 \sum_{n(q)} \sum_{x(q)*} e_q(nx+n\overline{x}+ln)=
 c\#\{
 x(c),\ x^2+lx+1\equiv 0 (c)\}.
\end{equation}
It is not equal to $R(l,n,c)$ in general. However the expression is very close. It is intriguing that for the evaluation of~\eqref{exp-MU} one sums over $n$ first while for the evaluation of~\eqref{def:R} one sums over $x$ first.

It is noticeable that \Lquote{one may view} (a shadow of) the Hecke relation:
\begin{equation}
 \lambda_f(p)^2=\lambda_f(p^2)-1
\end{equation}
at the level of exponential sums.
\subsection{Review of the geometric approach.}\label{sec:remarks:JK} In~\cites{Jorg01,JK06}, Jorgenson-Kramer developed a general approach via heat kernels that yields precise bounds for Arakelov invariants of Riemann surfaces. For the convenience of the reader we review briefly the main ingredients in their proof of a bound for $\CmZ_X$ when specialized to $N$ prime. A key difference in their approach and ours is that we apply Petersson trace formula and~\cites{Jorg01,JK06} apply Selberg trace formula. The case of a general level $N$ (square-free or not) is delicate and the reader is referred to~\cite{JK06} for a clear solution.

A lower bound for $\CmZ_{X}$ is obtained in~\cite{Jorg01}*{Theorem 3.3} as the logarithm of the sum of the volume, the number of cusps, the number of elliptic elements and the number of exceptional eigenvalues. As a consequence, in~\cite{Jorg01}*{\S 5.2} it is deduced a lower bound $\CmZ_{Y_0(N)}\ge O(\log N)$ where the constant in the $O()$ may be chosen roughly equal to $-4$.

An upper bound for $\CmZ_X$ is obtained in~\cite{Jorg01}*{Theorem 4.7},
see also the final formulas occurring in the proof the theorem and~\cite{Jorg01}*{Theorem 4.8}. It is given in terms of a short sum of inverse of exceptional eigenvalues, plus a weighted sum of length of short geodesics, plus a term controlling the remainder in the prime geodesic theorem (Huber constant). As a consequence, in~\cite{Jorg01}*{\S 5.3} it is deduced an upper bound $\CmZ_{Y_0(N)}\le O_\epsilon(N^\epsilon)$. It seems difficult to improve this growth from the geometric bound, in particular one needs both to control the number of short geodesics (which is done using results of Huxley) and the Huber constant (which is discussed \cite{Jorg01}*{Remark 4.10} and \cite{JK:geodesic}).

When $X$ is compact a nearby study is performed in~\cite{JK06}, in particular the following clean lower bound is proven: $\CmZ_X\ge -4\log (2g_X-2)$ where $g_X>1$ is the genus of $X$. An upper bound for $\CmZ_{X_0(N)}$ follows from \cite{JK06}*{Propositions 4.2, 5.3, 5.5}. In those propositions some of the steps consist of relating hyperbolic invariants of $X_0(N)$ to those of $Y_0(N)$. At present there is no direct relation between $\CmZ_{Y_0(N)}$ and the compactified $\CmZ_{X_0(N)}$. It would be interesting to investigate if the constructions given in~\cite{Broo99} may lead to such a relation.

\subsection{$\CmZ_X$ as a spectral invariant.}\label{sec:rem:spectral} The Selberg zeta function $Z(s)$ is a spectral invariant of the hyperbolic surface, see for instance Sarnak~\cite{Sarn87:det}. Following similar arguments, it is possible to state a precise identity between $\CmZ_X$ and the regularized trace of $\Delta^{-1}$, which we couldn't locate in the literature. We assume $X$ compact. Let $0=\lambda_0<\lambda_1\le \lambda_2 ...$ be the eigenvalues of $\Delta$ indexed in increasing order and counted with eventual multiplicity.

The standard way to regularize the trace is, see~\cite{Stei05}:
\begin{equation}\label{def-trace}
 \trd:=
 \lim_{w\rightarrow 1}
 \sum_{j\ge 1}
 \lambda^{-w}_j -\frac{1}{w-1}.
\end{equation}
The reader is referred to~\cite{Stei05} for the introduction of that trace and an interpretation as a \emph{geometrical mass}. Here $\{X_0(N)\}$ is a very particular case of Riemannian surfaces (arithmetic and hyperbolic). It is important to observe also that the genus tends to infinity with $N$. Usually one ought to study a fixed genus and let the metric varies, in which case the spectral zeta function may be used to measure the degeneration of the surface. Here since the genus is unbounded the significance of $\CmZ_X$ is not clear. Still for future convenience we believe it is useful to compute the difference $\MTr' \Delta-\CmZ_X$ explicitly.

There is a similarity between~\eqref{def-trace} and the regularization~\eqref{def:zx} that defines $\CmZ_X$. But since the spectral zeta function (which is $w\mapsto \MTr'(\Delta^{-w})$) and the Selberg zeta function (which is related to $w\mapsto \det(\Delta+s(s-1))$) are very different objects, the following identity is not obvious. We follow closely the arguments in~\cite{Sarn87:det} that provides a similar identity for the regularized \emph{determinant} of $\Delta$.
\begin{proposition} Let $X=\Gamma\SB \FmH$ be compact of genus $g$. Then:
\begin{equation}
 \trd=\CmZ_X +4(1+\gamma)(g-1)-1.
\end{equation}
Here $\gamma$ denotes the Euler constant.
\end{proposition}

\begin{proof} \footnote{To keep track the arguments it is worth to have in mind the following equality. If $\lambda>0$ and $f(w,s):=(\lambda+s(s-1))^{-w}$, then:
\begin{equation}
 -\frac{\partial f}{\partial w \partial s}(0,1)=\lambda^{-1}.
\end{equation}}
For $\MRe w >1$ and $s>1$, let (usual trace): 
\begin{equation}
 H(w,s):=\MTr(\Delta+s(s-1))^{-w}.
\end{equation}
It is known from the properties of the heat kernel that $H$ has a meromorphic continuation in $w\in \BmC$ and is regular in a neighborhood of $w=0$, see~\cite{Sarn87:det}*{(2.5)}. By definition~\eqref{def-trace}, we have:
\begin{equation}
 \trd =\lim_{w\to 1}
 \lim_{s\to 1}
 H(w,s)-(s(s-1))^{-w}-\frac{1}{w-1}.
\end{equation}
(observation: the two limits in $s$ and $w$ cannot be inverted at this stage).

The following identity is true for $\MRe w>1$ by absolute convergence and then for $\MRe w>0$ by holomorphic continuation (the left-hand-side is a priori meromorphic and the the right-hand side is already holomorphic):
\begin{equation}\label{H1+w}
 \frac{-1}{w(2s-1)}\frac{\partial}{\partial s}H(w,s)=H(1+w,s).
\end{equation}
Thus:
\begin{equation}\label{pf:tr}
 -\trd=\lim_{w\to 0}
 \lim_{s\to 1}
 \frac{1}{w(2s-1)}\frac{\partial}{\partial s}H(w,s)
+(s(s-1))^{-w-1}
+\frac{1}{w}.
\end{equation}

The spectral interpretation of the Selberg zeta function $Z(s)$ demands some care, particularly about the issue of constants. In the clear treatment given in~\cite{Sarn87:det}, the key input is that $\log Z(s)\to 0$ as $s\to \infty$. More precisely, differentiating~\cite{Sarn87:det}*{Theorem 1} \footnote{or from~\cite{Sarn87:det}*{(2.20)} together with~\cite{Sarn87:det}*{(2.15)}}, one has:
\begin{equation}
 \frac{Z'}{Z}(s)+\frac{Z'_\infty}{Z_\infty}(s)=\frac{d}{ds}\log \det (\Delta+s(s-1)) + F_0(2s-1),
\end{equation}
where by definition:
\begin{equation}
 \log \det (\Delta+s(s-1)) = -\frac{\partial H}{\partial w}(0,s).
\end{equation}
Whence:
\begin{equation}\label{pf:zx}
 \CmZ_X+\frac{Z'_\infty}{Z_\infty}(1)-F_0=
 -
 \lim_{s\to 1}
 \biggl(
 \frac{d}{ds} 
 \frac{\partial H}{\partial w}(0,s)
 +\frac{1}{s-1}
 \biggr).
\end{equation}
The constant $F_0$ is equal to $-(2g-2)$.

To conclude we would need to link the right-hand sides of~\eqref{pf:tr} and \eqref{pf:zx}. For this we use the following two facts:
\begin{equation}\tag{H1}
 H(w,s)-(s(s-1))^{-w}
\end{equation}
is regular for $(w,s)$ in a neighborhood of $(0,1)$.
\begin{equation}\tag{H2}
 \frac{d}{ds} H(0,s)=-1,
 \quad
 \forall s>1.
\end{equation}
The second fact (H2), because of~\eqref{H1+w}, is equivalent to the convergence of the limit in the definition of the regularized trace~\eqref{def-trace}. It is thus classical and follows for instance from~\cite{Sarn87:det}*{(2.1)}, see also~\cite{Stei05}*{appendix A}. The first fact (H1) follows from~\cite{Sarn87:det}*{(2.5) (2.6)}.

Now the proposition follows from~\eqref{pf:tr}, \eqref{pf:zx}, (H1) and (H2) by straightforward computations. Observe also that $2\frac{\Gamma'_2}{\Gamma_2}(1)=1-\log(2\pi)$ from~\cite{Sarn87:det}*{(1.13)} and $\frac{Z'_\infty}{Z_\infty}(1)=(2g-2)(1+2\gamma)$.
\end{proof}

\subsection{Squeezing the logarithms.}\label{sec:rem:log} In this paragraph we explain how to modify the arguments to get a control of logarithmic type.
\begin{theorem}\label{th:log}
 For $N$ prime, one has:
\begin{equation}
  \frac{1}{g(N)}\sum_{f\in \CmB_2(N)}
 \frac{L'}{L}(1,\MSym^2 f)
 =2\frac{\zeta'}{\zeta}(2)
 +O(\frac{(\log N)^7}{N}).
\end{equation}
\end{theorem}
\begin{proof} The scheme of the proof is the same as for the Theorem~\ref{th:main}.

In \S~\ref{sec:proof:poisson} we need to refine the dyadic subdivision. Instead of considering $1\le Y\le XN^\epsilon$, we force $1\le Y\le X$ and consider for the last function of the subdivision a function $U$ which is smooth of compact support on $(1,\infty)$ and is identically $1$ on $(2,\infty)$. Then the truncation is $U(\frac{n}{X})W_s(\frac{n}{X})$. This dyadic procedure costs a usual multiplicative factor $\log N$. A further $\log N$ is lost because of the asymptotic of $W_s(y)$ for $y$ small.

In the next step \S~\ref{sec:proof:truncation}, we simply cut $l$ up to $\abs{l}\le N^{20}=L$ with, say, two integrations by parts ($i=2$).

Then we integrate by parts only once ($i=1$) and insert the bound $J_1(x)\ll \min(1,x)$:
\begin{equation}
 g(l,c)\ll \frac{Y}{c}(1+\frac{Y}{c})
 \min(1,\frac{Y}{c})
 \frac{1}{1+\abs{l}}
 \ll \frac{Y}{x}
 \frac{1}{1+\abs{l}}.
\end{equation}

In \S~\ref{sec:proof:exp} we keep the factor $\tau(c)^2$ intact. Then the conclusion \S~\ref{sec:proof:conclusion} is modified into:
\begin{equation}
 B(X)\ll \frac{1}{X} \sum_{
 \substack{%
 c\equiv 0(N),\\ 1\le c\le C}}
 \sum_{\abs{l}\le L}
 \frac{X}{c}\frac{1}{1+\abs{l}} \tau(c)^2 \ll (\log N)^5
\end{equation}
we have used the fact that $\sum_{c\le C} \frac{\tau(c)^2}{c}\ll (\log C)^4$ which follows by a comparison of $\sum_{c\ge 1}\frac{\tau(c)^2}{c^s}$ with $\zeta(s)^4$ which has a pole of order $4$ at $s=1$. Collecting the logarithms we get the remainder term $O(\frac{(\log N)^7}{N})$ in Theorem~\ref{th:log}.
\end{proof}

\remark It is certainly possible to improve the exponent $7$ to some extent. Our aim was only to give to simple arguments yielding to a power of $\log N$. 

\subsection*{acknowledgements}
The author thanks the Institute for Advanced Study for its warm hospitality and Young-Heon Kim, Philippe Michel and Peter Sarnak for helpful comments.




\bibliographystyle{spbasic}      



\end{document}